\documentclass[12pt]{article}
\usepackage{amsmath}

\usepackage{graphicx}

\input{tcilatex}

\begin{document}

\title{COMPARISON OF DIFFERENT GOODNESS-OF-FIT TESTS}
\author{B. Aslan and G. Zech, Universit\"{a}t Siegen}
\maketitle

\begin{abstract}
Various distribution free goodness-of-fit test procedures have been
extracted from literature. We present two new binning free tests, the
univariate three-region-test and the multivariate energy test. The power of
the selected tests with respect to different slowly varying distortions of
experimental distributions are investigated. None of the tests is optimum
for all distortions. The energy test has high power in many applications and
is superior to the $\chi ^{2}$ test.
\end{abstract}

\section{INTRODUCTION}

Goodness-of-fit (gof) tests are designed to measure the compatibility of a
random sample with a theoretical probability distribution function (pdf).
The null hypothesis $H_{0}$ is that the sample follows the pdf. Under the
assumption that H$_{0}$ applies, the fraction of wrongly rejected
experiments - the probability of committing an error of the first kind - is
fixed to typically a few percent. A test is considered powerful if the
probability of accepting $H_{0}$ when $H_{0}$ is wrong - the probability of
committing an error of the second kind - is low. Of course, without
specifying the alternatives, the power cannot be quantified.

A discrepancy between a data sample and the theoretical description can be
of different origin. The problem may be in the theory which is wrong or the
sample may be biased by measurement errors or by background contamination.
In natural sciences we mainly have the latter situation. Even though the
statistical description is the same in both cases the choice of the specific
test may be different. In our applications we are mainly confronted with
``slowly varying'' deviations between data and theoretical description
whereas in other fields where for example time series are investigated,
``high frequency'' distortions are more likely.

Goodness-of-fit tests are based on classical statistical methods and are
closely related to classical interval estimation, but they contain also
Bayesian elements. Those, however, are only related with some prejudice on
the alternative hypothesis which affects the purity of the accepted
decisions and not the error of the first kind.

The power of one dimensional tests is not always invariant against
transformations of the variates. In more than one dimension (number of
variates), an invariant description is not possible.

Tests are classified in distribution dependent and distribution free tests.
The former are adapted to special pdfs like Gaussian, exponential or uniform
distributions. We will restrict our discussion mainly to distribution free
tests and tests which can be adapted to arbitrary distributions. Here we
distinguish tests applied to binned data and binning free tests. The latter
are in principle preferable but so far they are almost exclusively limited
to one dimensional distributions. A further distinction concerns the
alternative hypothesis. Usually, it is not restricted but there exist also
tests where it is parametrized.

Physicists tend to be content with $\chi ^{2}$ tests which are not
necessarily optimum in all cases. A very useful and comprehensive survey of
goodness-of-fit tests can be found in Ref. \cite{dago86} from 1986. Since
then, some new developments have occurred and the increase in computing
power has opened the possibility to apply more elaborate tests.

In Section 2 we summarize the most important tests. To keep this article
short we do not discuss \ tests based on the order statistic and spacing
tests. In Section 3 we introduce two new tests, the three region test and
the energy test. To compare the tests we apply them in Section 4 to some
specific alternative hypotheses. We do not consider explicitly composite
hypotheses.

\section{SOME RELEVANT TESTS}

\subsection{Chi-squared test}

The $\chi ^2$ test is closely connected to least square fits with the
difference that the hypothesis is fixed. The test statistic is 
\[
\chi ^2=\sum_{i=1}^B\frac{(Y_i-t_i)^2}{\delta _i^2}
\]
with the random variable $Y_i$, $t_i$ the expectation $E(Y_i)$ and $\delta
_i^2$ the expectation $E((Y_i-t_i)^2)$. Obviously the expectation value of $%
\chi ^2$ is 
\[
E(\chi ^2)=B
\]
In the Gaussian approximation $Y_i$ follows a Gaussian with mean $t_i$ and
variance $\delta _i^2$ and the test statistic follows a $\chi ^2$
distribution function $F_B(\chi ^2)$ with $B$ degrees of freedom. The
probability of an error of the first kind $\alpha $ (significance level,
p-value) defines $\chi _0^2$ with $F_B(\chi _0^2)=1-\alpha $. The null
hypothesis is rejected if in an actual experiment we find $\chi ^2>\chi _0^2$%
.

We obtain the Pearson test when the random variables $N_{i}$ are Poisson
distributed. 
\[
\chi ^{2}=\sum_{i=1}^{B}\frac{(N_{i}-t_{i})^{2}}{t_{i}}
\]

In the large sample limit the test statistic $\chi ^{2}$ has approximately a 
$\chi ^{2}$ distribution with $B$ degrees of freedom. When we have a
histogram, where a total of $N$ events are distributed according to a
multinomial distribution among the $B$ bins with probabilities $p_{i}$ we
get 
\[
\chi ^{2}=\sum_{i=1}^{B}\frac{(N_{i}-Np_{i})^{2}}{Np_{i}}
\]
which again follows asymptotically a $\chi ^{2}$ distribution, this time
with $B-1$ degrees of freedom. The reduced number of degrees of freedom is
due to the constraint $\Sigma N_{i}=N$.

Nowadays, the distribution function of the test statistic can be computed
numerically without much effort. The $\chi ^2$ test then can also be applied
to small samples. The Gaussian approximation is no longer required.

The $\chi ^{2}$ test is very simple and needs only limited computational
power. A big advantage compared to most of the other methods is that is can
be applied to multidimensional histograms. There are however also serious
drawbacks:

\begin{itemize}
\item  Its power in detecting slowly varying deviations of a histogram from
predictions is rather poor due to the neglect of possible correlations
between adjacent bins.

\item  Binning is required and the choice of the binning is arbitrary.

\item  When the statistics is low or the number of dimensions is high, the
event numbers per bin may be low. Then the asymptotic properties are no
longer valid and systematic deviations are hidden by statistical
fluctuations.
\end{itemize}

There are proposals to fix the bin widths by the requirement of equal number
of expected entries per bin. This is not necessarily the optimum choice \cite
{zechdesy}. Often there are outliers in regions where no events are expected
which would be hidden in wide bins.

For the number of bins a dependence on the sample size $n$%
\[
B=2n^{2/5}
\]
is proposed in Ref. \cite{moor86}. Our experience is that in most
experiments the number of bins is chosen too high. The sensitivity to slowly
varying deviations roughly goes with $B^{-1/4}$ \cite{zechdesy}. In
multidimensional cases the power of the test often can be increased by
applying it to the marginal distributions.

There is a whole class of $\chi ^{2}$ like tests. Many studies can be found
in the literature. The reader is referred to Ref. \cite{moor86}.

\subsection{Binning-free empirical distribution function tests}

The tests described in this section have been taken from the article by
Stephens in Ref. \cite{dago86}.

Supposing that a random sample of size $n$ is given, we form the order
statistic $X_{1}<X_{2}<...<X_{n}$. We consider the empirical distribution
function (EDF) 
\[
F_{n}(x)=\frac{\text{\# of observations }\leq x}{n}
\]
or

\begin{center}
$
\begin{array}{ll}
F_{n}(x)=0 & x<X_{1} \\ 
F_{n}(x)=\frac{i}{n} & X_{i}\leq x<X_{i+1} \\ 
F_{n}(x)=1 & X_{n}\leq x
\end{array}
$
\end{center}

$F_n(x)$ is a step function which is to be compared to the distribution $%
F(x) $ corresponding to $H_0$.

\begin{figure}
\centering
\includegraphics{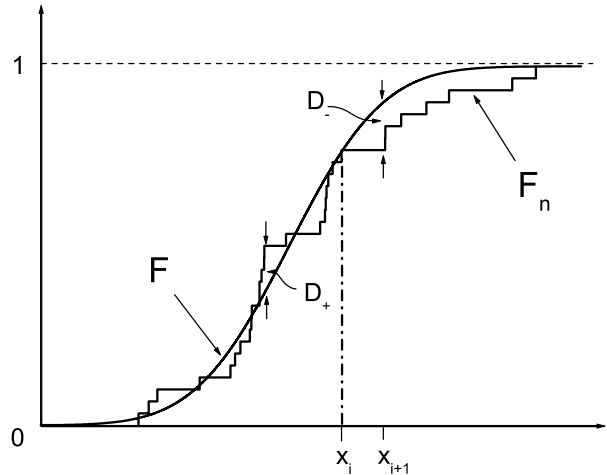}
\caption{Comparison of empirical and theoretical distributions}
\end{figure}

The EDF is consistent and unbiased. The tests discussed in this section are
invariant under transformation of the random variable. Because of this
feature, we can transform the distribution to the uniform distribution and
restrict our discussion to the latter.

\subsubsection{Probability integral transformation}

The probability integral transformation (PIT) 
\[
Z=F(X)
\]
transforms a general pdf $f(X)$ of $X$ into the uniform distribution $%
f^{*}(Z)$ of $Z$ . 
\begin{eqnarray*}
f^{*}(Z) &=&1;\text{ \ \ \ \ }0\leq Z\leq 1 \\
F^{*}(Z) &=&Z
\end{eqnarray*}
The underlying idea of this transformation is that the new EDF of $Z$, $%
F^{*}(Z)$ is extremely simple and that it conserves the distribution of the
test quantities discussed in this section. It is easily seen that

\[
F_n(x)-F(x)=F_n^{*}(z)-z
\]

Note, however, that the PIT does not necessarily conserve all interesting
features of the gof problem. Resolution effects are washed out and for
example in a lifetime distribution, an excess of events at small and large
lifetimes may be judged differently but are treated similarly after a PIT.
It is not logical to select specific gofs for specific applications but to
transform all kinds of pdfs to the same uniform distribution. The PIT is
very useful because it permits standardization but one has to be aware of
its limitations.

\subsubsection{Supremum statistics}

The maximum positive (negative) deviation of $F_n(x)$ from $F(x)$ $D_{+}$ ($%
D_{-}$) (see Fig. 1) are used as tests statistics. Kolmogorov
(Kolmogorov-Smirnov test) has proposed to use the maximum absolute
difference. Kuiper uses the sum $V=D_{+}+D_{-}$. This test statistic is
useful for observations ``on the circle'' for example for azimuthal
distributions where the zero angle is a matter of definition. 
\begin{eqnarray*}
D_{+} &=&\sup_x\left\{ F_n(x)-F(x)\right\} \\
D_{-} &=&\sup_x\left\{ F(x)-F_n(x)\right\} \\
D &=&\sup_x\left\{ \left| F_n(x)-F(x)\right| \right\} \ \ \ \ \ \ \ \ \ \ \ 
\text{(Kolmogorov)} \\
V &=&D_{+}+D_{-}\ \ \ \ \ \ \ \ \ \ \ \ \ \ \ \ \ \ \ \ \ \ \ \ \ \ \ \ \ \ 
\text{(Kuiper)}
\end{eqnarray*}

The supremum statistics are invariant under the PIT.

\subsubsection{Integrated deviations - quadratic statistics}

The Cramer-von Mises family of tests measures the integrated quadratic
deviation of $F_{n}(x)$ from $F(x)$ suitably weighted by a weighting
function $\psi $: 
\[
Q=n\int_{-\infty }^{\infty }\left[ F(x)-F_{n}(x)\right] ^{2}\psi (x)dF
\]
In the standardized form we have 
\begin{equation}
Q=n\int_{0}^{1}\left[ z-F_{n}^{\ast }(z)\right] ^{2}\psi (z)dz  \label{qdef}
\end{equation}

Since the construction of $F_{n}^{\ast }(Z)$ includes already an
integration, $F_{n}^{\ast }(z_{i})$ and $F_{n}^{\ast }(z_{k})$ are not
independent and the additional integration in Equation \ref{qdef} is not
obvious.

With $\psi _{CvM}=1$ we get the Cramer-von Mises statistic $W^{2}$ and $\psi
_{AD}=\left[ z(1-z)\right] ^{-1}$ leads to the Anderson-Darling statistic $%
A^{2}$. 
\begin{eqnarray*}
W^{2} &=&n\int_{0}^{1}\left[ z-F_{n}^{\ast }(z)\right] ^{2}dz\ \ \ \ \ \ \ \
\ \ \ \text{(Cramer-von Mises)} \\
A^{2} &=&n\int_{0}^{1}\frac{\left[ z-F_{n}^{\ast }(z)\right] ^{2}}{z(1-z)}%
dz\ \ \ \ \ \ \ \ \ \ \ \text{(Anderson - Darling)}
\end{eqnarray*}
The Anderson-Darling statistic $A^{2}$ weights strongly deviations near $z=0$
and $z=1$. This is justified because there the experimental deviations are
small due to the constraints $\left[ z-F_{n}^{\ast }(z)\right] =0$ at $z=0$
and $z=1$.

Watson has proposed a quadratic statistic on the circle: 
\[
U^{2}=n\int_{0}^{1}\left\{ F_{n}^{\ast }(z)-z-\int_{0}^{1}\left[ F_{n}^{\ast
}(z)-z\right] dz\right\} ^{2}dz\ \ \ \ \ \ \ \ \ \ \ \text{(Watson)}
\]

\subsection{The Neyman statistic test}

This test is different from all previously discussed tests. It parametrizes
the alternative hypothesis and applies the likelihood ratio test. The
alternative hypothesis corresponds to a pdf of the exponential family: 
\[
g_{k}(z)=C(\theta _{1},\theta _{2}...\theta _{k})\exp \left[
\sum_{i=1}^{k}\theta _{i}\pi _{i}(z)\right]
\]
$g_{k}(z)$ are smooth alternatives to uniformity. The functions $\pi _{i}$
are Legendre polynomials of order $i$, $\theta _{i}$ are free parameters and 
$C$ is a normalization function. The number $k$ of parameters is selected by
the user.

The likelihood ratio leads to the test statistic 
\[
N_k=\frac 1n\sum_{i=1}^k\left( \sum_{j=1}^n\pi _i(z_j)\right) ^2
\]

Asymptotically, for large values of $N_{k}$, $N_{k}$ is distributed
according to the $\chi ^{2}$ distribution with $k$ degrees of freedom.

\section{NEW TESTS}

\subsection{Three region test}

Often experimental distributions are biased by an excess or lack of events
in certain regions of the random variable. We have designed a test which
subdivides the variable space into three pieces, containing $%
n_{1},n_{2},n_{3}=n-n_{1}-n_{2}$ events, such that the deviation between
data and prediction from $H_{0}$ is maximum. \ The test quantity is 
\[
O=\sup_{n_{1},n_{2}}\left\{
w_{1}(n_{1}-np_{1})^{2}+w_{2}(n_{2}-np_{2})^{2}+w_{3}(n_{3}-np_{3})^{2}%
\right\}
\]
where $np_{k}$ are the expectation values and $w_{k}$ weights depending on $%
np_{k}$. The specific choice 
\begin{eqnarray*}
w_{k} &=&\frac{1}{np_{k}} \\
O_{\chi } &=&\sup_{n_{1},n_{2}}\left\{ \frac{(n_{1}-np_{1})^{2}}{np_{1}}+%
\frac{(n_{2}-np_{2})^{2}}{np_{2}}+\frac{(n_{3}-np_{3})^{2}}{np_{3}}\right\}
\end{eqnarray*}
maximizes $\chi ^{2}$ of the three bins. In the comparison below we have
chosen weights equal to one.

Of course the test can be generalized to a higher number of subregions.

\subsection{Minimum energy test}

\subsubsection{The idea}

Let us assume that we have a continuous charge distribution $\rho (\vec{r})$
of positive electric charges and a sample of negative point charges with
total charge equal to minus the integrated positive charge. The potential
energy is minimum when the negative point charges follow $\rho $. Then, up
to effects due to the discrete nature of the point charges, the charge
density is zero everywhere. Any displacement of charges would increase the
energy. We use this property to construct a binning free test procedure.

We simulate the theoretical distribution by $m$ charges of charge $1/m$
each. Usually, these charges are distributed using a Monte Carlo simulation.
To the $n$ experimental sample points (data points) we associate charges $%
-1/n$. The test quantity $\phi $ corresponds to the potential energy. It
contains two terms $\phi _{1},\phi _{2}$ corresponding to the interaction of
the experimental charges with each other and to the interaction of the
experimental charges with the positive simulation charges. 
\begin{align}
\phi _{1}& =\frac{1}{n^{2}}\sum_{i<j}R(d_{ij})  \label{phi1} \\
\phi _{2}& =-\frac{1}{nm}\sum_{i,j}R(t_{ij})  \label{phi2} \\
\phi & =\phi _{1}+\phi _{2}  \label{phi}
\end{align}

Here $d_{ij}$ is the distance between two data points and $t_{ij}$ is the
distance between a data point and a simulation point and $R$ is a
correlation function defined below. The sums run over all combinations.

\textbf{Remark}: The minimum energy requirement for the equality of
experimental and theoretical distribution is strictly correct only when the
number $m$ of simulation charges is equal to the number $n$ of experimental
charges. For the general case with a continuous theoretical distribution or
simulation sample and experimental sample of different size, the optimum
agreement of the two distributions is not well defined and there is a slight
dependence of the minimum energy configuration on the correlation function.
This is however a purely academic problem, the test statistic $\phi$ remains
a powerful indicator for an incompatibility of the experimental sample with $%
H_{0}$.

\subsubsection{The correlation function}

We note that the minimum energy configuration does not depend on the
application of the one-over-distance power law of electrostatics. We may
apply a wide class of correlation functions $R(r)$ with the only requirement
that $R$ has to decrease monotonically with the distance $r$.

We have investigated three different types of \ correlation functions, power
laws, a logarithmic dependence and Gaussians. 
\begin{align}
R_{\kappa}(r) & =\frac{1}{r^{\kappa}}  \label{rpot} \\
R_{l}(r) & =-\ln r  \label{rlog} \\
R_{s}(r) & =e^{-r^{2}/(2s^{2})}  \label{rexp}
\end{align}

The first type is motivated by the analogy to electrostatics, the second is
long range and the third emphasizes a limited range for the correlation
between different points. The power $\kappa$ of the denominator in Equ. \ref
{rpot} and the parameter $s$ in Equ. \ref{rexp} may be chosen differently
for different dimensions of the sample space and different applications. For
long range distortions a small value of $\kappa$ around $0.1$ is
recommended. For short range deviations the test quantity with larger values
around $0.3$ is more sensitive.

The inverse power law and the logarithm have a singularity at $r$ equal to
zero. Very small distances, however, should not be weighted too strongly
since distortions with sharp peaks are not expected and usually inhibited by
the finite experimental resolution. We eliminate the singularity by
introducing a lower cutoff $d_{\min }$ for the distances $d$ and $t$.
Distances less than $d_{\min }$ are replaced by $d_{\min }$. The value of
this parameter is not critical, it should be of the order of the average
distance $d$ in the regions where the $f_{0}$ is maximum and not less than
the experimental resolution.

The energy test with Gaussian correlation function is closely related to the
Pearson $\chi ^{2}$ test. A more detailed description of the energy test is
presented in Ref. \cite{asla02}.

\subsection{Comparison of uni-variate tests}

\begin{figure}
\centering
\includegraphics{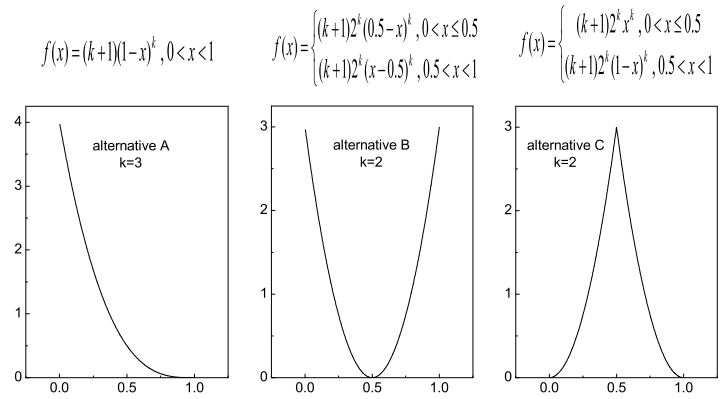}
\caption{\label{background1d} Different types of background distributions}
\end{figure}

We have tested the null hypothesis of a uniform distribution in the interval 
$[0,1]$ using a uniform distribution contaminated by the background
distributions displayed in Figure \ref{background1d}.

Background hypothesis A modifies the mean, hypotheses B, C change the
variance of $H_{0}.$

\begin{figure}
\centering
\includegraphics{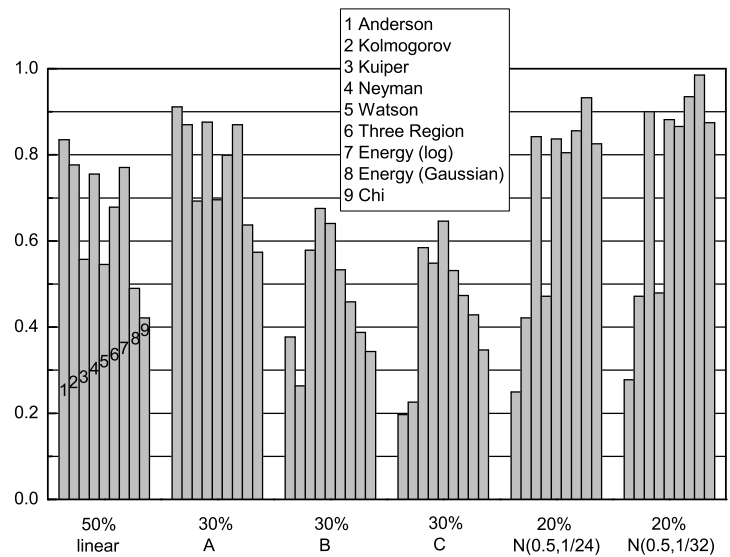}
\caption{\label{power1d} Histograms show the powers of different tests at 5\% level for the sample size n=100.}
\end{figure}

The power of various tests described above is presented in Figure \ref{power1d}.

As expected, non of the tests is optimum for all kind of distortions.
Several tests perform better than the $\chi ^{2}$ test. The Neyman test, the
Anderson-Darling test and the Kolmogorov-Smirnov test are sensitive to a
shift in the mean. Anderson's test detects especially deviations at the
borders of the interval. Watson's and Kuiper's tests are useful for the
detection of distortions of the variance. The two new tests compare
favorably with the standard ones.

\section{MULTIVARIATE TESTS}

The Mardia test \cite{mardia} and the Neyman smooth test \cite{rayn89} can
be used to investigate two-dimensional Gaussian distributions. The only
distribution free test known to us which is independent of the dimensions of
the variate space is the $\chi ^{2}$ test. A generalized Kolmogorov-Smirnov
test \cite{kolmogmulti} depends on the ordering scheme of the variates. The
binning free energy test developed by us is also independent of the number
of variates, however the distribution of test statistic has to be computed
for the specific sample distribution under study.

\subsection{Comparison of multivariate tests}

\begin{figure}
\centering
\includegraphics{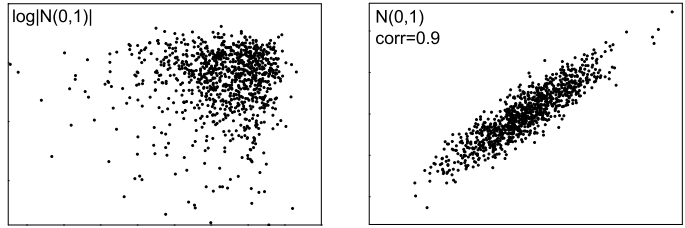}
\caption{\label{background2d} Different types of background distributions in two dimensions}
\end{figure}

We have used a two-dimensional Gaussian null hypothesis and contaminated the
sample with the background distributions shown in Figure \ref{background2d}.

\begin{figure}
\centering
\includegraphics{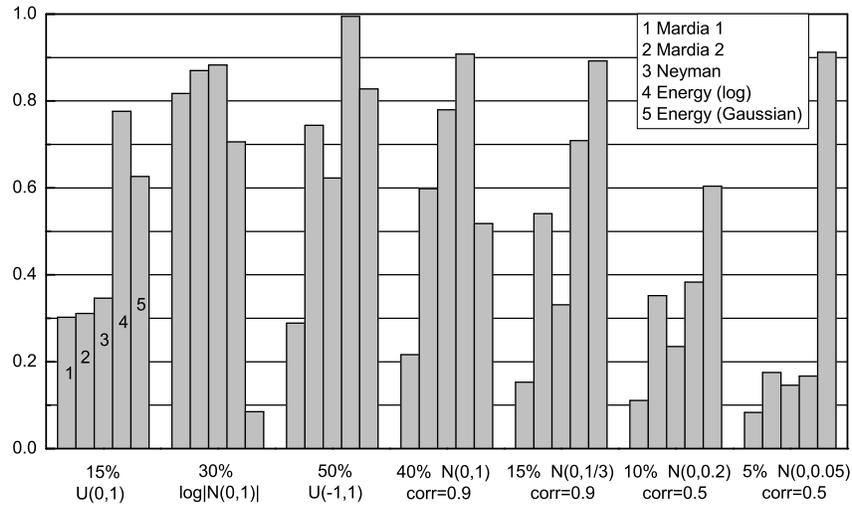}
\caption{\label{power2d} Power comparisons of tests of bivariate normality for the sample size n=200 at 5\% significance level. }
\end{figure}

The power of the two Mardia tests, the Neyman smooth test and the energy
test with logarithmic and Gaussian correlation function is presented in
Figure \ref{power2d}.

In most cases the two energy tests perform better than the alternatives even
though those have been designed for a Gaussian null hypothesis.

\subsection{Example: Comparing experimental data to a Monte Carlo prediction}

\begin{figure}
\centering
\includegraphics{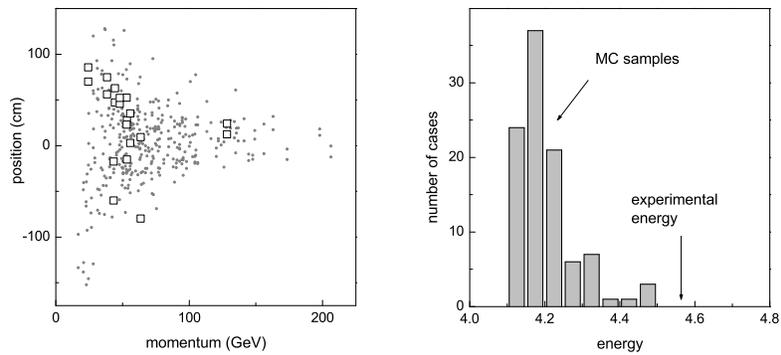}
\caption{\label{scatter} Comparison of experimental
distribution (squares) with Monte Carlo simulation (dots). The experimental
energy computed from the scatter plot (left) is compared to a Monte Carlo
simulation of the experiment (right).}
\end{figure}

\allowbreak In Figure \ref{scatter}, left hand side, we compare the position
and momentum of a few $J/\psi $ decay tracks to a Monte Carlo simulation.The
right hand plot compares the energy computed from the distribution on the
left hand side to a Monte Carlo simulation of the null hypothesis. The
experimental point, indicated by the arrow, is larger than all Monte Carlo
values. Apparently, the data do not follow the prediction.

\section{ \ CONCLUSIONS}

The $\chi ^{2}$ test suffers from the requirement to choose a binning. In
one dimension it should be replaced by the well established binning free
tests like the Kolmogorov test. The choice of a specific test has to depend
on the expected kind of possible distortion of the theoretical distribution.
For a localized background we advise to use the new three region test. For
multivariate applications the new energy test is an attractive alternative
to the $\chi ^{2}$ test.


\begin{thebibliography}{9}
\bibitem{dago86}  ``Goodness-of-fit techniques'', ed. R. B. D'Agostino and
M. A. Stephens, Dekker (1986).

\bibitem{zechdesy}  G. Zech, ``Comparing statistical data to Monte Carlo
simulation - parameter fitting and unfolding'', Desy 95-113 (1995).

\bibitem{moor86}  D. S. Moore, ``Tests of Chi Squared Type'', in Ref. 1, p
63.

\bibitem{asla02}  B. Aslan and G. Zech, ``A new class of binning free,
multivariate goodness-of-fit tests: the energy tests'', hep-ex/0203010
(2002).

\bibitem{mardia}  V. K. Mardia, ``Measures of multivariate skewness and
kurtosis with applications'', Biometrika 57 (1970) 519.

\bibitem{rayn89}  J. C. W. Rayner et al, ``Smooth Tests of Goodness of
Fit'', Oxford University Press, Oxford (1989).

\bibitem{kolmogmulti}  A. Justel et al, ``A multivariate Kolmogorov-Smirnov
test of goodness of fit'', Statist. \& Prob. Lett. 35 (1997) 251.
\end{thebibliography}
\end{document}